\documentclass[a4paper, 10pt, leqno]{article}
\usepackage{amsmath, amssymb}

\newtheorem{theorem}{\bf Theorem}
\newtheorem{proposition}[theorem]{\bf Proposition}
\newtheorem{lemma}[theorem]{\bf Lemma}
\newtheorem{corollary}[theorem]{\bf Corollary}

\def\epsilon{\varepsilon}

\begin{document}

%\noindent Running title: Expanding fronts in an anisotropic diffusion equation  

%\vspace*{2.4em} 

%\vspace*{1.2em}
%\newpage

\Large \noindent 
{\bf Remarks on space-time behavior in the Cauchy problems  
of the heat equation and the curvature flow equation 
with mildly oscillating initial values}  

\vspace*{0.8em}

\normalsize
\noindent Hiroki Yagisita%${}^*$

\vspace*{2.4em}

\noindent {\bf Abstract} \ 
We study two initial value problems of the linear diffusion equation $u_t=u_{xx}$ 
and the nonlinear diffusion equation $u_t=(1+{u_x}^2)^{-1}u_{xx}$, 
when Cauchy data $u(x,0)=u_0(x)$ are bounded and oscillate mildly. 
The latter nonlinear heat equation is the equation of the curvature flow,  
when the moving curves are represented by graphs. 
In the case of $\lim_{|x|\rightarrow+\infty} |x{u_0}^\prime(x)|=0$, 
by using an elementary scaling technique, we show 
$$\lim_{t\rightarrow+\infty}\left|u(\sqrt t x,t)
-\left(F(-x)u_0(-\sqrt t)+F(+x)u_0(+\sqrt t)\right)\right|=0$$
for the linear heat equation $u_t=u_{xx}$, 
where $x\in\mathbb R$ and $F(z):=\frac{1}{2\sqrt \pi}\int_{-\infty}^z e^{-\frac{y^2}{4}} dy$. 
Further, by combining with a theorem of Nara and Taniguchi, 
we have the same result for the curvature equation $u_t=(1+{u_x}^2)^{-1}u_{xx}$.  
In the case of $\lim_{|x|\rightarrow +0} |x{u_0}^\prime(x)|=0$ 
and in the case of $\sup_{x\in\mathbb R} |x{u_0}^\prime(x)|<+\infty$, respectively, 
we also give a similar remark for the linear heat equation $u_t=u_{xx}$. 

\vspace*{0.8em} 


\noindent nontrivial dynamics, 
nontrivial large-time behavior, irregular behavior.

\vfill

%\noindent 
%A proposed running title: Traveling Waves in Nonlocal Systems II. 

\vspace*{1.6em}

\footnoterule

%\noindent 
%Communicated by H. Okamoto. 
%Received November 1, 2008. 
%Revised February 28, 2009. 

%\noindent 
%2000 Mathematics Subject Classifications: 35K57, 35K65, 35K90, 45J05.

\noindent
%${}^*$ 
\ Department of Mathematics, 
Faculty of Science, 
Kyoto Sangyo University

%\noindent
%\ Motoyama, Kamigamo, Kita-Ku, Kyoto-City, 603-8555, Japan 

\newpage

\section{Introduction} 
\( \, \, \, \, \, \, \, \) 
In this paper, by using an elementary scaling argument, 
we study space-time behavior 
in the Cauchy problem of the heat equation  
\begin{equation}
\left\{\begin{array}{l}
u_t(x,t)=u_{xx}(x,t), \ \ \ (x,t) \in \mathbb R \times (0,+\infty), \\ 
u(x,0)=u_0(x), \ \ \ x \in \mathbb R, 
\end{array}\right.
\end{equation}
when the initial values $u_0(x)$ are bounded and oscillate mildly. 
We also study the Cauchy problem of the nonlinear diffusion equation 
\begin{equation}
\left\{\begin{array}{l}
u_t(x,t)=\frac{u_{xx}(x,t)}{1+(u_x(x,t))^2}, \ \ \ (x,t) \in \mathbb R \times (0,+\infty), \\ 
u(x,0)=u_0(x), \ \ \ x \in \mathbb R, 
\end{array}\right. 
\end{equation} 
which is the equation of the curvature flow when the moving curves are represented by graphs.

First, we mention criteria for stabilization 
of the solution $u(x,t)$ to the Cauchy problem of the heat equation $u_t=u_{xx}$. 
From [3, 11, 4, 2] (e.g.), we see the following: 
\begin{theorem}
Let $u_0\in L^\infty(\mathbb R)$ and $c\in\mathbb R$. Then, 
the solution $u(x,t)$ to (1.1) satisfies $\lim_{t\rightarrow+\infty}u(x,t)=c$ 
if and only if $u_0(x)$ satisfies $\lim_{R\rightarrow+\infty}\frac{1}{2R}\int_{-R}^{+R}u_0(x+y)dy=c$. 
Moreover, $u(x,t)$ satisfies $\lim_{t\rightarrow+\infty}\sup_{x\in\mathbb R}|u(x,t)-c|=0$ 
if and only if $u_0(x)$ 
satisfies $\lim_{R\rightarrow+\infty}\sup_{x\in\mathbb R}|\frac{1}{2R}\int_{-R}^{+R}u_0(x+y)dy-c|=0$. 
\end{theorem}
On the other hand, Collet and Eckmann [1] given a simple example of a bounded initial value $u_0(x)$ where the solution $u(x,t)$ to (1.1) 
oscillates forever as $t\rightarrow+\infty$ 
(See also [5]). 
So, the large-time behavior of a solution $u(x,t)$ to (1.1) 
with a bounded initial value $u_0(x)$ may be complex. 
Indeed, V\'{a}zquez and Zuazua [13] showed the general behavior is very complex: 
\begin{theorem} \ 
(i) \ Let $u_0\in L^\infty(\mathbb R)$. 
Then, the set of accumulation points in $L^\infty_{loc}(\mathbb R)$ 
of $\{(e^{\triangle t} u_0)(\sqrt t \, \cdot \, )\}_{t>0}$ as $t\rightarrow+\infty$ 
coincides with the set 
$\{(e^\triangle\phi)( \, \cdot \, ) \, | \, \phi\in A\}$, 
where $A$ is the set 
of accumulation points of $\{u_0(\lambda \, \cdot \, )\}_{\lambda>0}$ 
as $\lambda\rightarrow+\infty$ in the weak-star topology $\sigma(L^\infty,L^1)$. 

(ii) \ Let $c>0$ and $B_c=\{f\in L^\infty(\mathbb R) \, | \, \|f\|_{L^\infty}\leq c\}$. 
Let $\mathcal M_c$ be the set of $f\in B_c$ 
such that the set of accumulation points  
of $\{f(\lambda \, \cdot \, )\}_{\lambda>0}$ as $\lambda\rightarrow+\infty$
in the weak-star topology $\sigma(L^\infty,L^1)$ 
is $B_c$. Then, 
$\mathcal M_c$ is dense with empty interior in $B_c$ 
with the weak-star topology $\sigma(L^\infty,L^1)$. 
\end{theorem}
They also showed the general behavior in a number of evolution equations 
on $\mathbb R^N$ is complex. However, the behavior may be rather simple, 
if the initial value oscillates mildly.  
In this paper, we prove the following, which is a remark 
on the long-time behavior in the Cauchy problem (1.1) 
when the initial value $u_0(x)\in L^\infty(\mathbb R)\cap C^1(\mathbb R\setminus\{0\})$ 
satisfies $\lim_{|x|\rightarrow+\infty}|x{u_0}^\prime(x)|=0$: 
\begin{theorem}
Let $u_0\in L^\infty(\mathbb R)\cap C^1(\mathbb R\setminus\{0\})$ 
and $\lim_{|x|\rightarrow+\infty} |x{u_0}^\prime(x)|=0$.  
Then, the solution $u(x,t)$ to (1.1) satisfies  
$$\lim_{t\rightarrow+\infty}\sup_{x\in [-L,+L]}\left|u(\sqrt t x,t)
-\left(F(-x)u_0(-\sqrt t)+F(+x)u_0(+\sqrt t)\right)\right|=0$$
for all $L>0$, where $F(z):=\frac{1}{2\sqrt \pi}\int_{-\infty}^z e^{-\frac{y^2}{4}} dy$. 
\end{theorem}
\begin{corollary}
Let $u_0\in L^\infty(\mathbb R)\cap C^1(\mathbb R\setminus\{0\})$ 
and $\lim_{|x|\rightarrow+\infty} |x{u_0}^\prime(x)|=0$.  
Then, the set of accumulation points in $L^\infty_{loc}(\mathbb R)$ 
of $\{(e^{\triangle t} u_0)(\sqrt t \, \cdot \, )\}_{t>0}$ 
as $t\rightarrow+\infty$ 
coincides with the set 
$\{\alpha F(- \, \cdot \, )+\beta F(+ \, \cdot \, ) \, | \, (\alpha,\beta)\in A\}$, 
where $F(z):=\frac{1}{2\sqrt \pi}\int_{-\infty}^z e^{-\frac{y^2}{4}} dy$ 
and $A$ is the set of accumulation points in ${\mathbb R}^2$ of $\{(u_0(-\lambda),$ 
%%%%%%%%%%%%%%%%%%%%%%%%
%%%%%%%%%%%%%%%%%%%%%%%%
%%%%%%%%%%%%%%%%%%%%%%%%
$u_0(+\lambda))\}_{\lambda>0}$ 
as $\lambda\rightarrow+\infty$. 
\end{corollary}
We also prove the following two: 
\begin{proposition}
Let $u_0\in L^\infty(\mathbb R)\cap C^1(\mathbb R\setminus\{0\})$ 
and $\lim_{|x|\rightarrow+0} |x{u_0}^\prime(x)|=0$.  
Then, the solution $u(x,t)$ to (1.1) satisfies  
$$\lim_{t\rightarrow+0}\sup_{x\in [-L,+L]}\left|u(\sqrt t x,t)
-\left(F(-x)u_0(-\sqrt t)+F(+x)u_0(+\sqrt t)\right)\right|=0$$
for all $L>0$, where $F(z):=\frac{1}{2\sqrt \pi}\int_{-\infty}^z e^{-\frac{y^2}{4}} dy$. 
\end{proposition}
\begin{proposition}
Let $u_0\in C^1(\mathbb R\setminus\{0\})$ 
and $\sup_{x\in\mathbb R\setminus\{0\}} |x{u_0}^\prime(x)| < +\infty$.  
Then, the solution $u(x,t)$ to (1.1) satisfies  
$$\left|u(\sqrt t x,t)
-\left(F(-x)u_0(-\sqrt t)+F(+x)u_0(+\sqrt t)\right)\right|$$
$$\leq G(-x)\left(\sup_{y<0} |y{u_0}^\prime(y)|\right) 
+ G(+x)\left(\sup_{y>0} |y{u_0}^\prime(y)|\right)$$
for all $(x,t)\in\mathbb R \times (0,+\infty)$, where 
$F(z):=\frac{1}{2\sqrt \pi}\int_{-\infty}^z e^{-\frac{y^2}{4}} dy$ 
and $G(z):=\frac{1}{2\sqrt \pi}\int_{0}^{+\infty} 
e^{-\frac{(z-y)^2}{4}} |\log y|dy$. 
\end{proposition}

\noindent 
{\bf Remark} \ {\rm (i)} \ 
Let $(a, b) \in {\mathbb R}^2$ and 
$$u_0(x)=\left\{ \begin{array}{ll}
a \ & (x<0), \\ 
b \ & (x>0).
\end{array} \right.$$
Then, the solution $u(x,t)$ to (1.1) satisfies  
$$u(\sqrt t x,t) = a F(-x) + b F(+x)$$ 
for all $(x,t)\in\mathbb R \times (0,+\infty)$, where 
$F(z):=\frac{1}{2\sqrt \pi}\int_{-\infty}^z e^{-\frac{y^2}{4}} dy$. 

{\rm (ii)} \ Let $u(x,t)$ be the solution to (1.1). 
Then, the function $$v(x,t):=u(e^{\frac{t}{2}}x, e^t)$$ 
is the solution to 
$$
\left\{\begin{array}{l}
v_t(x,t)=v_{xx}(x,t)+\frac{x}{2}v_{x}(x,t), 
\ \ \ (x,t) \in {\mathbb R}^2, \\ 
v(x,0)=(e^{\triangle} u_0)(x), \ \ \ x \in \mathbb R. 
\end{array}\right.
$$ 

{\rm (iii)} \ Let $u_1(x)=\phi_1(\log (-x))$ and $u_2(x)=\phi_2(\log (+x))$. 
Then, $x {u_1}^\prime(x)={\phi_1}^\prime (\log (-x))$ 
and $x {u_2}^\prime(x)={\phi_2}^\prime (\log (+x))$. 

\vspace*{1.2em}

\noindent 
Nara and Taniguchi [9] showed that the difference between the solution to the heat equation (1.1) 
and that to the curvature flow equation (1.2) with the same initial value 
is of order $O(t^{-\frac{1}{2}})$ as $t\rightarrow+\infty$. 
Precisely, they given the following theorem: 
\begin{theorem}
Let $\varepsilon>0$. Suppose $u_0\in C^2(\mathbb R)$ satisfies 
$\sup_{x\in\mathbb R}(|u_0(x)|+|{u_0}^\prime(x)|+|{u_0}^{\prime\prime}(x)|)<+\infty$ 
and $\sup_{x_1,x_2\in\mathbb R, x_1\not=x_2}
\frac{|{u_0}^{\prime\prime}(x_1)-{u_0}^{\prime\prime}(x_2)|}{|x_1-x_2|^\varepsilon}<+\infty$. 
Then, the maximum interval of existence of the classical solution $u(x,t)$ to (1.2) 
is $[0,+\infty)$ and the solution $u(x,t)$ satisfies 
$$\sup_{t>0,x\in\mathbb R}t^\frac{1}{2}\left|u(x,t)-
\frac{1}{2\sqrt{\pi t}}\int_{-\infty}^{+\infty}e^{-\frac{(x-y)^2}{4t}}u_0(y)dy\right|<+\infty.$$
\end{theorem}
Therefore, by combining it with Theorem 3, we have the following remark 
on the long-time behavior in the Cauchy problem (1.2):   
\begin{corollary}
Let $\varepsilon>0$. Suppose $u_0\in C^2(\mathbb R)$ satisfies 
$\sup_{x\in\mathbb R}(|u_0(x)|+|{u_0}^\prime(x)|$ 
%%%%%%%%%%%%%%%%%%%%%%%%%%%%%%%%%%
%%%%%%%%%%%%%%%%%%%%%%%%%%%%%%%%%%
%%%%%%%%%%%%%%%%%%%%%%%%%%%%%%%%%%
$+|{u_0}^{\prime\prime}(x)|)<+\infty$, 
$\sup_{x_1,x_2\in\mathbb R, x_1\not=x_2}
\frac{|{u_0}^{\prime\prime}(x_1)-{u_0}^{\prime\prime}(x_2)|}{|x_1-x_2|^\varepsilon}<+\infty$ 
and $\lim_{|x|\rightarrow+\infty}$ 
%%%%%%%%%%%%%%%%%%%%%%%%%%%%%%%%%%%%%%%%%
%%%%%%%%%%%%%%%%%%%%%%%%%%%%%%%%%%%%%%%%%
%%%%%%%%%%%%%%%%%%%%%%%%%%%%%%%%%%%%%%%%%
$|x{u_0}^\prime(x)|=0$. 
Then, the solution $u(x,t)$ to (1.2) 
satisfies 
$$\lim_{t\rightarrow+\infty}\sup_{x\in [-L,+L]}\left|u(\sqrt t x,t)
-\left(F(-x)u_0(-\sqrt t)+F(+x)u_0(+\sqrt t)\right)\right|=0$$
for all $L>0$, where $F(z):=\frac{1}{2\sqrt \pi}\int_{-\infty}^z e^{-\frac{y^2}{4}} dy$. 
\end{corollary}

Nara [8] showed that the difference between the solution to the heat equation on $\mathbb R^N$ 
and that to the mean curvature flow equation on $\mathbb R^N$ with the same initial value 
is of order $O(t^{-\frac{1}{2}})$ as $t\rightarrow+\infty$, 
when the initial value is radially symmetric. 
See [12, 6] for the difference 
between the behavior of a disturbed planar front 
in a bistable reaction-diffusion equation 
and that of a mean curvature flow with a drift term. 
See [10, 14, 12, 7] for other nontrivial large-time behaviors 
in nonlinear diffusion equations. 

\newpage

\section{Proof}
\begin{lemma} 
The solution $u(x,t)$ to (1.1) satisfies 
$$\sup_{x\in[-L,+L]}\left|u(\sqrt t x,t)
-\left(aF(-x)+bF(+x)\right)\right|$$
$$\leq \frac{1}{2\sqrt{\pi}} \int_{0}^{+\infty} 
\rho_{L}(z) 
\left(\left|u_0(-\sqrt t z)-a\right| 
+ \left|u_0(+\sqrt t z)-b\right|\right) 
dz$$ 
for all $(L,t)\in (0,+\infty)^2$ and $(a,b)\in {\mathbb R}^2$, 
where $\rho_{L}(z):=\sup_{z_0\in[-L,+L]}e^{-\frac{(z-z_0)^2}{4}}$. 
\end{lemma} 
{\bf Proof.} \ From 
$$u(\sqrt t x,t)=\frac{1}{2\sqrt{\pi}}\int_{-\infty}^{+\infty} 
e^{-\frac{(x-y)^2}{4}}u_0(\sqrt{t}y) dy$$
$$=\frac{1}{2\sqrt{\pi}}\int_{0}^{+\infty} 
e^{-\frac{(x+z)^2}{4}}u_0(-\sqrt{t}z) dz
+\frac{1}{2\sqrt{\pi}}\int_{0}^{+\infty} 
e^{-\frac{(x-z)^2}{4}}u_0(+\sqrt{t}z) dz,$$ 
we see 
$$u(\sqrt t x,t)
-\left(aF(-x)+bF(+x)\right)
$$
$$= \frac{1}{2\sqrt{\pi}}\int_{0}^{+\infty} 
e^{-\frac{(z-(-x))^2}{4}}\left(u_0(-\sqrt{t}z)-a\right)dz$$ 
$$+ \frac{1}{2\sqrt{\pi}}\int_{0}^{+\infty} 
e^{-\frac{(z-(+x))^2}{4}}\left(u_0(+\sqrt{t}z)-b\right)dz.$$
So, we have the conclusion. 
\hfill 
$\blacksquare$

\begin{lemma} 
Let $u_0\in C^1(\mathbb R\setminus\{0\})$ and $\alpha>0$. 
Then, $\lim_{|x|\rightarrow+\infty} |x{u_0}^\prime(x)|=0$ 
implies $\lim_{|s|\rightarrow+\infty}|u_0(s\alpha)-u_0(s)|=0$. 
Also, $\lim_{|x|\rightarrow+0} |x{u_0}^\prime(x)|=0$ 
implies $\lim_{|s|\rightarrow+0}|u_0(s\alpha)-u_0(s)|=0$. 
\end{lemma} 
{\bf Proof.} \ We see
$$|u_0(s\alpha)-u_0(s)|
=\left|\int_{1}^{\alpha} s{u_0}^\prime(sz) dz\right|$$ 
$$\leq \left(\alpha+\frac{1}{\alpha}\right)
\sup_{\min\{\alpha,\frac{1}{\alpha}\} \leq |z| 
\leq \max\{\alpha,\frac{1}{\alpha}\}} 
|s{u_0}^\prime(sz)|$$ 
$$\leq \left(\alpha+\frac{1}{\alpha}\right)^2 
\sup_{\min\{\alpha,\frac{1}{\alpha}\} \leq |z| 
\leq \max\{\alpha,\frac{1}{\alpha}\}} |sz{u_0}^\prime(sz)|$$ 
$$= \left(\alpha+\frac{1}{\alpha}\right)^2
\sup_{\min\{\alpha,\frac{1}{\alpha}\}|s| \leq |x| 
\leq \max\{\alpha,\frac{1}{\alpha}\}|s|} |x{u_0}^\prime(x)|.$$ 
So, we have the conclusion. 
\hfill 
$\blacksquare$ 

\vspace*{0.8em}

\noindent 
{\bf Proof of Theorem 3 and Proposition 5.} \  We see 
$$\left|u_0(-\sqrt t z)-u_0(-\sqrt t)\right| 
+ \left|u_0(+\sqrt t z)-u_0(+\sqrt t)\right| 
\leq 4\|u_0\|_{L^\infty(\mathbb R)}$$ 
for all $t>0$ and $z>0$. Hence, because of 
$\rho_{L} \in L^1({\mathbb (0,+\infty)})$, 
we have the conclusions by Lemmas 9 and 10. 
\hfill 
$\blacksquare$ 

\vspace*{0.8em}

\noindent 
{\bf Remark} \ 
Let $u_1(x)=\phi_1(\log (-x))$, $u_2(x)=\phi_2(\log (+x))$ 
and $\alpha>0$. 
Then, $\lim_{|z|\rightarrow+\infty} |\phi_1(z+\log \alpha)-\phi_1(z)|=0$ 
implies $\lim_{s\rightarrow-\infty}|u_1(s\alpha)-u_1(s)|=0$ 
and $\lim_{s\rightarrow-0} |u_1(s\alpha)-u_1(s)|=0$. 
Also, $\lim_{|z|\rightarrow+\infty} |\phi_2(z+\log \alpha)-\phi_2(z)|=0$ 
implies $\lim_{s\rightarrow+\infty}|u_2(s\alpha)-u_2(s)|=0$ 
and $\lim_{s\rightarrow+0} |u_2(s\alpha)-u_2(s)|=0$.

\vspace*{0.8em}

\noindent 
{\bf Proof of Proposition 6.} \ From 
$$u(\sqrt t x,t)
-\left(F(-x)u_0(-\sqrt t)+F(+x)u_0(+\sqrt t)\right)
$$
$$= \frac{1}{2\sqrt{\pi}}\int_{0}^{+\infty} 
e^{-\frac{(x+z)^2}{4}}\left(u_0(-\sqrt{t}z)-u_0(-\sqrt{t})\right)dz$$
$$+ \frac{1}{2\sqrt{\pi}}\int_{0}^{+\infty} 
e^{-\frac{(x-z)^2}{4}}\left(u_0(+\sqrt{t}z)-u_0(+\sqrt{t})\right)dz$$ 
$$= \frac{1}{2\sqrt{\pi}}\int_{0}^{+\infty} 
e^{-\frac{(x+z)^2}{4}} \left( 
\int_{1}^{z} \frac{(-\sqrt{t}y){u_0}^\prime(-\sqrt{t}y)}{y} dy 
\right)dz$$
$$+ \frac{1}{2\sqrt{\pi}}\int_{0}^{+\infty} 
e^{-\frac{(x-z)^2}{4}}\left( 
\int_{1}^{z} \frac{(+\sqrt{t}y){u_0}^\prime(+\sqrt{t}y)}{y} dy 
\right)dz,$$ 
we see 
$$\left|u(\sqrt t x,t)
-\left(F(-x)u_0(-\sqrt t)+F(+x)u_0(+\sqrt t)\right)\right|
$$
$$\leq \frac{1}{2\sqrt{\pi}}\int_{0}^{+\infty} 
e^{-\frac{(x+z)^2}{4}} \left| 
\int_{1}^{z} \frac{|(-\sqrt{t}y){u_0}^\prime(-\sqrt{t}y)|}{y} dy 
\right|dz$$
$$+ \frac{1}{2\sqrt{\pi}}\int_{0}^{+\infty} 
e^{-\frac{(x-z)^2}{4}}\left| 
\int_{1}^{z} \frac{|(+\sqrt{t}y){u_0}^\prime(+\sqrt{t}y)|}{y} dy 
\right|dz$$ 
$$\leq \frac{1}{2\sqrt{\pi}}\int_{0}^{+\infty} 
e^{-\frac{(x+z)^2}{4}} \left| 
\int_{1}^{z} \frac{\sup_{s<0}|s{u_0}^\prime(s)|}{y} dy 
\right|dz$$
$$+ \frac{1}{2\sqrt{\pi}}\int_{0}^{+\infty} 
e^{-\frac{(x-z)^2}{4}}\left| 
\int_{1}^{z} \frac{\sup_{s>0}|s{u_0}^\prime(s)|}{y} dy 
\right|dz$$ 
$$= \frac{\sup_{s<0}|s{u_0}^\prime(s)|}{2\sqrt{\pi}}\int_{0}^{+\infty} 
e^{-\frac{((-x)-z)^2}{4}} |\log z| dz$$
$$+ \frac{\sup_{s>0}|s{u_0}^\prime(s)|}{2\sqrt{\pi}}\int_{0}^{+\infty} 
e^{-\frac{((+x)-z)^2}{4}} |\log z| dz.$$ 
So, we have the conclusion. 
\hfill 
$\blacksquare$ 

%\newpage 

\vspace*{1.6em} 
%%%%%%%%%%%%%%%%%%%%%
%%%%%%%%%%%%%%%%%

\noindent Acknowledgments. \ 
It was partially supported 
by Grant-in-Aid for Scientific Research (No. 23540254) 
from Japan Society for the Promotion of Science. 

\newpage 

%%%%%%%%%%%%%%%%%%%%%%%%%%%
%%%%%%%%%%%%%%%%%%
%\vspace*{0.8em} 
%%%%%%%%%%%%%%%%%%%%%%%%%%%
%%%%%%%%%%%%%%%%%%
%%%%%%%%%%%%%%%%%%%%%%%%%%%
%%%%%%%%%%%%%%%%%%

\[ \begin{array}{c} \mbox{R\scriptsize EFERENCES}  \end{array} \]
 
[1] P. Collet and J.-P. Eckmann, 
Space-time behaviour in problems of hydrodynamic type: a case study, 
{\it Nonlinearity}, 5 (1992), 1265-1302. 

[2] V. N. Denisov and V. D. Repnikov, The stabilization 
of a solution of a Cauchy problem for parabolic equations, 
{\it Differential Equations}, 20 (1984), 16-33. 

[3] S. D. Eidel'man, {\it Parabolic systems}, North-Holland, 1969. 

[4] S. Kamin, On stabilisation of solutions of the Cauchy problem 
for parabolic equations, {\it Proc. Roy. Soc. Edinburgh A}, 
76 (1976), 43-53. 

[5] M. Krzy\.{z}a\'{n}ski, Sur l'allure asymptotique des potentiels 
de chaleur et de l'int\'{e}grale de Fourier-Poisson, 
{\it Ann. Polon. Math.}, 3 (1957), 288-299. 

[6] H. Matano and M. Nara, Large time behavior of disturbed planar fronts 
in the Allen-Cahn equation, {\it J. Differential Equations}, 
251 (2011), 3522-3557. 

[7] H. Matano, M. Nara and M. Taniguchi, Stability of planar waves 
in the Allen-Cahn Equation, {\it Comm. Partial Differential Equations}, 
34 (2009), 976-1002. 

[8] M. Nara, Large time behavior of radially symmetric surfaces 
in the mean curvature flow, {\it SIAM J. Math. Anal.}, 
39 (2008), 1978-1995. 

[9] M. Nara and M. Taniguchi, The condition on the stability 
of stationary lines in a curvature flow in the whole plane, 
{\it J. Differential Equations}, 237 (2007), 61-76. 

[10] P. Pol\'{a}\v{c}ik and E. Yanagida, Nonstabilizing solutions 
and grow-up set for a supercritical semilinear diffusion equation, 
{\it Differential Integral Equations}, 17 (2004), 535-548. 

[11] V. D. Repnikov and S. D. \`{E}\u{i}del'man, 
A new proof of the theorem on the stabilization of the solution 
of the Cauchy problem for the heat equation, 
{\it Math. USSR-Sbornik}, 2 (1967), 135-139. 

[12] J.-M. Roquejoffre and V. R.-Michon, Nontrivial large-time behaviour 
in bistable reaction-diffusion equations, {\it Annali di Matematica}, 
188 (2009), 207-233. 

[13] J. L. V\'{a}zquez and E. Zuazua, 
Complexity of large time behaviour of evolution equations with bounded data, 
{\it Chin. Ann. of Math. B}, 23 (2002), 293-310. 

[14] E. Yanagida, Irregular behavior of solutions for Fisher's equation, 
{\it J. Dyn. Diff. Equat.}, 19 (2007), 895-914.  

%\newpage

%602-0904

%‹ž"sŽsã‹ž‹æ¬"‡'¬ 538-301

%–ö‰º_‹I

\end{document}